\documentclass{amsart}

\usepackage{amsmath}
\usepackage{amsthm}
\usepackage{hyperref}
\usepackage{amsfonts,graphics,amsthm,amsfonts,amscd,latexsym}
\usepackage{epsfig}
\usepackage{flafter}
\usepackage{mathtools}
\usepackage{comment}
\usepackage{stmaryrd}

\usepackage{mathabx,epsfig}

\hypersetup{
    colorlinks=true,    
    linkcolor=blue,          
    citecolor=blue,      
    filecolor=blue,      
    urlcolor=blue           
}
\usepackage{tikz}
\usetikzlibrary{graphs,positioning,arrows,shapes.misc,decorations.pathmorphing}

\tikzset{
    >=stealth,
    every picture/.style={thick},
    graphs/every graph/.style={empty nodes},
}

\tikzstyle{vertex}=[
    draw,
    circle,
    fill=black,
    inner sep=1pt,
    minimum width=5pt,
]
\usepackage[position=top]{subfig}
\usepackage{amssymb}
\usepackage{color}

\calclayout

\usetikzlibrary{decorations.pathmorphing}
\tikzstyle{printersafe}=[decoration={snake,amplitude=0pt}]

\newcommand{\rank}{\operatorname{rank}}

\newcommand{\PP}{\mathbb{P}}

\newcommand{\CC}{\mathbb{C}}
\newcommand{\GG}{\mathbb{G}}
\newcommand{\xX}{{\mathcal X}}
\newcommand{\yY}{{\mathcal Y}}
\newcommand{\dD}{{\mathcal D}}
\newcommand{\kK}{{\mathcal K}}
\newcommand{\OO}{{\mathcal O}}
\newcommand{\EE}{{\mathcal E}}
\newcommand{\ver}{\operatorname{vert}}
\newcommand{\sym}{\mathop{\mathrm{Sym}}\nolimits}

\def\O#1.{\mathcal {O}_{#1}}			
\def\pr #1.{\mathbb P^{#1}}				
\def\af #1.{\mathbb A^{#1}}			
\def\ses#1.#2.#3.{0\to #1\to #2\to #3 \to 0}	
\def\xrar#1.{\xrightarrow{#1}}			
\def\K#1.{K_{#1}}						
\def\bA#1.{\mathbf{A}_{#1}}			
\def\bM#1.{\mathbf{M}_{#1}}				
\def\bL#1.{\mathbf{L}_{#1}}				
\def\bB#1.{\mathbf{B}_{#1}}				
\def\bK#1.{\mathbf{K}_{#1}}			
\def\subs#1.{_{#1}}					
\def\sups#1.{^{#1}}

\usepackage{tikz}
\usetikzlibrary{matrix,arrows,decorations.pathmorphing}

  \newtheorem{introthm}{Theorem}

  \newtheorem{theorem}{Theorem}[section]

  \newtheorem{conjecture}[theorem]{Conjecture}
  \newtheorem{setup}{Setup}

    \theoremstyle{definition}
  \newtheorem{lemma}[theorem]{Lemma}
  \newtheorem{proposition}[theorem]{Proposition}

  \newtheorem{definition}[theorem]{Definition}

\theoremstyle{remark}

\numberwithin{equation}{section}

\usepackage[all]{xy}

\begin{document}

\title[The algebraic Green-Griffiths-Lang conjecture]{The algebraic Green-Griffiths-Lang conjecture for complements of very general pairs of divisors}

\author[K.~Ascher]{Kenneth Ascher}
\address[Ascher]{Department of Mathematics, University of California, Irvine, CA 92697, USA
}
\email{kascher@uci.edu}

\author[A.~Turchet]{Amos Turchet}
\address[Turchet]{Dipartimento di Matematica e Fisica, Universit\'a Roma Tre, I-00146 Roma, Italy
}
\email{amos.turchet@uniroma3.it}

\author[W. Yeong]{Wern Yeong}
\address[Yeong]{Department of Mathematics, University of California, Los Angeles, CA 90095, USA
}
\email{wyyeong@math.ucla.edu}

\begin{abstract}
 We prove that the complement of a very general pair of hypersurfaces of total degree $2n$ in $\PP^n$ is algebraically hyperbolic modulo a proper closed subvariety. 
This provides evidence towards conjectures of Lang--Vojta and Green--Griffiths, 
and partially extends previous work of Chen, Pacienza--Rousseau, and Chen--Riedl and the third author.
\end{abstract}

\subjclass[2020]{14J70, 14C20, 32Q45}
\keywords{Algebraic hyperbolicity, log variety, Brody hyperbolicity}

\maketitle

\section{Introduction}

The motivation of this paper originates from the study of hyperbolic properties of varieties of log general type. 
A major open problem in the area is the logarithmic Green--Griffiths--Lang conjecture,  which posits that if $(X,D)$ is a pair of log general type, then $X \setminus D$ is \emph{pseudo Brody hyperbolic}, i.e. there is a proper closed subset $Z \subsetneq X$ containing the images of all nonconstant holomorphic maps $\mathbb{C} \to X \setminus D$. 
Demailly in \cite{Demailly} introduced an algebraic analogue of Brody hyperbolicity in the compact case (i.e. when $D=\emptyset$), which was naturally extended to the following definition in the log case by Chen \cite{Chen}. 

\begin{definition} 
{ 
Let $X$ be a smooth projective variety and $D$ an snc divisor. 
We say that the pair $(X,D)$ is \emph{algebraically hyperbolic} if there exists an ample line bundle $L$ and a constant $\varepsilon > 0$ such that for every reduced, irreducible curve $Y \subset X$ such that $Y \not\subset D$, one has
\begin{equation}\label{eqn:ah-definition}
2g(Y)-2+\left|\nu_Y^{-1}(Y \cap D)\right| \geq \varepsilon \deg \nu_Y^* L,
\end{equation}
where $\nu_Y: Y^\nu \to Y$ is the normalization and $g(Y) := g(Y^\nu)$ denotes the geometric genus of $Y$. 

Let $Z \subsetneq X$ be a proper subvariety not contained in $D$.
We say that $(X,D)$ is \emph{algebraically hyperbolic modulo $Z$} if \eqref{eqn:ah-definition} holds for all curves $Y \not\subset D \cup Z$.
We call such a proper subvariety $Z$ an \textit{exceptional set} for the pair $(X,D)$.
}
\end{definition}

Brody hyperbolicity implies algebraic hyperbolicity for compact varieties and 
these hyperbolicity properties are conjectured to be equivalent \cite{Demailly}. 
In this vein, it is then natural to consider the following.

\begin{conjecture}[Algebraic Green-Griffiths-Lang Conjecture]\label{conj:LVFF}
{ \em 
A pair $(X,D)$ of log general type is algebraically hyperbolic modulo a proper, closed subvariety $Z \subsetneq X$.
}
\end{conjecture}

This conjecture is sometimes referred to as the \emph{Lang--Vojta conjecture for function fields} (see \cite{CZ08}), since it can be related to a special case of Vojta's Conjecture \cite[Conjecture 24.3]{Vojta} using the analogy between function fields and number fields (see for example \cite[Section 1.1]{CT} for an extended discussion).  

The case $X = \PP^n$ has attracted a lot of interest, being one of the first natural cases to consider. Historically the conjecture has been studied mostly either when $D$ has many irreducible components, or when $D$ is taken to be very general in the space of hypersurfaces of $\PP^n$. 
In the former case, Conjecture~\ref{conj:LVFF} is proven in \cite{CZGm} for $X = \PP^2$ and $D$ a curve with at least three irreducible components (see also \cite{CZ08,T,CT}). 
In these cases, the exceptional set $Z$ has been described in \cite{CT24}. 
In higher dimensions, \cite{Wang_etal} proved the conjecture for $X= \PP^n$ and $D$ a hypersurface of $\deg(D) \geq n+2$ with $n+1$ irreducible components, and more generally for any ramified cover of $\GG_m^n$. The corresponding hyperbolicity statements were previously obtained in \cite{NWY} proving the logarithmic Green-Griffiths-Lang Conjecture in these setting.

In the latter case, \cite{Chen} and \cite{PacienzaRousseau2007} proved independently that the complement of a \textit{very general} hypersurface of degree at least $2n+1$ in $\PP^n$ is algebraically hyperbolic, and in particular Conjecture~\ref{conj:LVFF} holds in these cases with $Z = \emptyset$. 
More recently, \cite{CRY} proved Conjecture~\ref{conj:LVFF} in the case where $D$ is a very general hypersurface in $X = \mathbb{P}^n$ of degree $2n$, and identified the exceptional locus $Z$ in the case $n = 2$. 

Our main result deals with the case when $D$ is very general in the space of hypersurfaces with two irreducible components.
\begin{introthm}\label{introthm}
{\em
Let $D\subset \PP^n$ be 
a very general pair of hypersurfaces of degrees $d_1,d_2$ such that $d_1+d_2=2n$. 
Let $H$ denote the hyperplane bundle on $\mathbb{P}^n$. Then there exists a proper subvariety $Z\subsetneq \PP^n$ such that 
every nonconstant map $f\colon Y \rightarrow \PP^n$ from a smooth projective curve $Y$ such that $f(Y)\nsubset Z\cup D$ satisfies
\begin{enumerate}
\item If $n \geq 3$, $2g(Y)-2+\left|f^{-1}(D)\right| \geq \deg f^*H$.
\item If $n=2$, $2g(Y)-2+\left|f^{-1}(D)\right| \geq \frac{1}{2}\deg f^*H$.
\end{enumerate}
}
\end{introthm}

We expect that the methods used to prove Theorem \ref{introthm} can be adapted to the case of arbitrarily many components ($n\geq 3$) or to prove log algebraic hyperbolicity of two or more components of total degree at least $2n+1$.

\subsection{Sketch of proof}
We follow the strategy of proof in \cite{CRY}, adapted to the case when $D$ consists of two components.
If, for a general pair $D$ of hypersurfaces of degrees $d_1$ and $d_2$, there is a map $f:Y\rightarrow \PP^n$ as above from a smooth projective curve $Y$ of geometric genus $g$ and $\deg f^*H=e$ satisfying $2g-2+\left|f^{-1}(D)\right|<e$ (or $2g-2+\left|f^{-1}(D)\right|<\frac{1}{2}e$ when $n=2$), then we may ``spread them out'' into appropriate families of curves and hypersurfaces using a standard construction (\S\ref{sec:hilb-scheme}).
We then study the log tangent sheaves of these families and relate them to Lazarsfeld kernel bundles (\S\ref{sec:kernel-bundles}). 
In \S\ref{sec:3}, assuming $d_1+d_2=2n$, we identify a special family of lines in $\PP^n$ passing through general points on such curves, and use them to prove that such curves form a bounded family which cannot dominate $\PP^n$.

\subsection*{Acknowledgements}
The authors would like to thank Angelo Felice Lopez, Emanuela Marangone 
and Eric Riedl
for many helpful discussions. KA was supported in part by NSF grants DMS-2140781 and DMS-2302550, and a UCI Chancellor's Fellowship.  AT is partially supported by PRIN  2022HPSNCR: Semiabelian varieties, Galois representations and related Diophantine problems and   PRIN   2020KKWT53: Curves, Ricci flat Varieties and their Interactions, and is a member of the INDAM group GNSAGA.

\section{Preliminaries}

In this section, after reviewing the notions of log tangent and log normal sheaves, we show that the bound appearing in Theorem~\ref{introthm} is naturally related to the degree of the log normal sheaf (Proposition~\ref{prop:degree-log-normal}).
We then recall the construction from \cite{PacienzaRousseau2007, CRY} of ``spreading out'' a general curve and hypersurface arrangement into appropriate families of curves and hypersurfaces. 
Finally, we prove some new results in \S\ref{sec:kernel-bundles} which will allow us to analyze these families via \textit{Lazarsfeld kernel bundles} (see Definition~\ref{def:Lazarsfeld}) in a similar manner to \cite{CRY}.

\subsection{Log tangent sheaves}\label{sec:log-tangent-sheaves}  
We briefly recall the construction and some properties of log tangent sheaves. 
We refer the reader to \cite{Saito1980, SernesiBook} for a more comprehensive treatment. 
Let $D\subset \PP^n$ be a degree $d$ reduced divisor. 
We define the {\em log tangent sheaf $T_{\PP^n}(-\log D)$ of $D$ in $\PP^n$} as the kernel of the natural map $T_{\PP^n}\rightarrow O_D(D).$ 
It sits in the four-term exact sequence
\begin{equation}\label{eqn:4term}
0 \rightarrow T_{\PP^n}(-\log D) \rightarrow T_{\PP^n} \rightarrow \OO_D(D) \rightarrow \OO_{D_\text{sing}}(D) \rightarrow 0,
\end{equation}
where $D_{\text{sing}}$ is the subscheme in $\PP^n$ defined by the equation for $D$ and its partials, and $\OO_{D_\text{sing}}(D)$ is called the \emph{first cotangent sheaf supported on $D_{\text{sing}}$}. 
When $D$ is the union of two plane curves in general position, $D_{\text{sing}}$ is the scheme consisting of the reduced points where the curves intersect.
Let us denote the kernel of the restriction map $\OO_D(D)\rightarrow \OO_{D_\text{sing}}(D)$ by $N'_{D/\PP^n}$, which is called the \emph{equisingular normal sheaf of $D$ in $\PP^n$}. 
Then \eqref{eqn:4term} shortens to the short exact sequence
\begin{equation}\label{eqn:log-tgt-fiber-1}
    0 \rightarrow T_{\PP^n}(-\log D) \rightarrow T_{\PP^n} \rightarrow N'_{D/\PP^n} \rightarrow 0.
\end{equation}
The log tangent sheaf of a reduced plane curve is always a vector bundle. 
In general, the log tangent sheaf of a reduced divisor in a smooth projective variety is only a reflexive sheaf. 
It is a vector bundle when the divisor is smooth or snc. 

Let $f: Y\rightarrow \PP^n$ be a nonconstant map from a smooth projective curve $Y$ such that $f(Y)\not\subset D$. 
Let $f^{-1}(D)\subset Y$ denote the reduced subscheme. 
As above, we may define the log tangent sheaf $T_Y(-\log D)$ of $f^{-1}(D)$ in $Y$, which sits inside the short exact sequence
\begin{equation}\label{eqn:log-tgt-fiber-2}
    0 \rightarrow T_Y(-\log D) \rightarrow T_Y \rightarrow \OO_{f^{-1}(D)}(f^{-1}(D))\rightarrow 0.
\end{equation}
Finally, we define the \emph{log normal sheaf $N_{f/\PP^n}(\log D)$ of $f$ with respect to $D$} as the quotient sheaf in the following short exact sequence:
\begin{equation}\label{eqn:log-normal-fiber}
    0 \rightarrow T_Y(-\log D) \rightarrow f^*T_{\PP^n}(-\log D)\rightarrow N_{f/\PP^n}(\log D)\rightarrow 0.
\end{equation}
The short exact sequences \eqref{eqn:log-tgt-fiber-1}, \eqref{eqn:log-tgt-fiber-2} and \eqref{eqn:log-normal-fiber} fit into the exact, commuting diagram
\begin{equation}\label{eqn:log-normal-diag}
    \xymatrix{
     & &0 \ar[d] &0 \ar[d] &0 \ar[d] & \\
        &0 \ar[r] &T_Y(-\log D) \ar[d] \ar[r] &T_Y \ar[r] \ar[d] &\OO_{f^{-1}(D)}(f^{-1}(D)) \ar[r] \ar[d] &0 \\
        &0 \ar[r] &f^*T_{\PP^n}(-\log D) \ar[d] \ar[r] &f^*T_{\PP^n} \ar[d] \ar[r] &f^*N'_{D/\PP^n} \ar[d] \ar[r] &0\\
        &0 \ar[r] &N_{f/\PP^n}(\log D) \ar[d] \ar[r] &N_{f/\PP^n} \ar[d] \ar[r] &K \ar[d] \ar[r] &0\\		
        & &0 &0 &0 &
    }
\end{equation}
where the middle column is the usual short exact sequence associated to $f: Y\rightarrow \PP^n$, and $K$ is the quotient sheaf in the right column (and bottom row). 

The following proposition shows that the bound in Theorem~\ref{introthm} is intimately related to the degree of the log normal sheaf $N_{f/\PP^n}(\log D)$.

\begin{proposition}[cf. Proposition 3.2 in \cite{CRY}]
\label{prop:degree-log-normal}
Let $D\subset \PP^n$ be a degree $d$ reduced divisor, and let $f\colon Y\rightarrow \PP^n$ be a nonconstant map from a smooth projective curve $Y$ such that $f(Y)\not\subset D$. 
Then
\begin{equation*}
    \deg N_{f/\PP^n}(\log D)= 2g-2+\left|f^{-1}(D)\right|-(d-n-1)\deg f^*H,
\end{equation*}
where $g=g(Y)$ is the geometric genus of $Y$ and $H$ is the hyperplane divisor in $\PP^n.$
\end{proposition}

\begin{proof}
We compute the degree of $N_{f/\PP^n}(\log D)$ using the diagram~\eqref{eqn:log-normal-diag}. 
The top row shows that \linebreak $\deg T_Y(-\log D) = 2-2g-\left|f^{-1}(D)\right|$. 
By definition, we have
\begin{equation*}
0 \rightarrow N'_{D/\PP^n} \rightarrow \OO_D(D) \rightarrow \OO_{D_\text{sing}}(D) \rightarrow 0,
\end{equation*}
so $\deg f^*N'_{D/\PP^n} = \deg f^*\OO_D(D) - \deg f^*\OO_{D_\text{sing}}(D) = d \cdot \deg f^*H.$
The middle row in the diagram then shows that $\deg f^*T_{\PP^n}(-\log D) = -(d-n-1) \cdot \deg f^*H.$
Finally, the result follows from taking degrees in the left column of the diagram.
\end{proof}

\subsection{Construction of families}\label{sec:hilb-scheme}
Suppose that for a divisor $D\subset \PP^n$ consisting of a general pair of (smooth and irreducible) hypersurfaces of degrees $d_1$ and $d_2$, there is a nonconstant map $f:Y \rightarrow \PP^n$ from a smooth projective curve $Y$ of geometric genus $g=g(Y)$ such that 
$f(Y)\nsubset D$, $\left|f^{-1}(D)\right|=i$, and $\deg f^*H=e$ for some $i,e.$ 
As in \cite{Voisin1996, ClemensRan2004, PacienzaRousseau2007, CRY}, we construct some families in which the above arrangement is a general fiber.

Let $B_1=\bigoplus\limits_{i=1}^2 H^0(\PP^n,\OO_{\PP^n}(d_i))$. 
By passing to an \'etale cover $B_2$ of $B_1$, we may first construct the families $\xX_2,\dD_2$ and $\yY_2$ over the base $B_2$, described as follows:
\begin{enumerate}
\item $\xX_2=\PP^n\times B_2$.
\item $\dD_2\subset \xX_2$ is the universal hypersurface over $B_2$.
\item $\yY_2\subset \xX_2$ is a family over $B_2$ of pointed curves of geometric genus $g$ and degree $e$ such that the number of intersection points between a general fiber of $\yY_2$ and of $\dD_2$ is $i$.
\end{enumerate}
After taking a resolution of $\yY_2/B_2$ and possibly shrinking $B_2$, we obtain families $\xX , \dD$ and $\yY$ over the base $B$, such that 
\begin{itemize}
\item $\xX = \PP^n \times B$,
\item $\dD$ is a snc divisor in $\xX$, 
\item the reduced divisor $f^{-1}(\dD)$ in $\yY$ is snc, and 
\item $\yY$ is a smooth family with a generically injective morphism $f: \yY \rightarrow \xX$ such that $f^{-1}(\dD)$ is a disjoint union of sections of $\yY/B$.
\end{itemize}
Let us denote by $p_1: \xX \rightarrow B$ and $p_2: \xX \rightarrow \PP^n$ the natural projection maps from $\xX$, and by $g: \dD \rightarrow \xX$ and $h:f^{-1}(\dD) \rightarrow \yY$ the natural injective morphisms over $B$.
We denote the fibers of $\dD$ and $\yY$ over $b\in B$ by $D_b$ and $Y_b$, respectively, 
the restriction of $f$ by $f_b:Y_b\rightarrow \PP^n$, and 
the inclusion of the fiber $Y_b$ in $\yY$ by $\iota_b:Y_b\hookrightarrow \yY$.

As in \S\ref{sec:log-tangent-sheaves}, we may define the log tangent sheaves $T_\xX(-\log \dD)$ and $T_\yY(-\log \dD)$ as follows:
\begin{equation*}\label{eqn:log-tgt-family-1}
    0 \rightarrow T_{\xX}(-\log \dD) \rightarrow T_{\xX} \rightarrow N'_{\dD/\xX} \rightarrow 0
\end{equation*}
\begin{equation*}\label{eqn:log-tgt-family-2}
    0 \rightarrow T_{\yY}(-\log \dD) \rightarrow T_{\yY} \rightarrow \OO_{f^{-1}(\dD)}(f^{-1}(\dD)) \rightarrow 0.
\end{equation*}
Since $(\xX,\dD)$ and $(\yY, f^{-1}(\dD))$ are snc pairs, the sheaves $T_\xX(-\log \dD)$ and $T_\yY(-\log \dD)$ are vector bundles. 
As before, we may define the log normal sheaf $N_{\yY/\xX}(\log \dD)$ on $\yY$ as the quotient sheaf in the following short exact sequence:
\begin{equation*}\label{eqn:log-normal-family}
    0 \rightarrow T_{\yY}(-\log \dD) \rightarrow f^*T_{\xX}(-\log \dD) \rightarrow N_{\yY/\xX}(\log \dD) \rightarrow 0.
\end{equation*}
Without loss of generality, we may assume that the families $\dD$ and $\yY$ are stable under the $GL(n+1)$-action, so that we have surjective maps $T_\xX(-\log \dD) \rightarrow p_2^* T_{\PP^n}$ and $T_\yY(-\log \dD) \rightarrow f^*p_2^*T_{\PP^n}$. 
We define the {\em vertical log tangent sheaves} $T_\xX^{\ver}(-\log \dD)$ and $T_\yY^{\ver}(-\log \dD)$ as the kernels in the following short exact sequences:
\begin{equation*}\label{eqn:verticalx}
0 \rightarrow T_\xX^{\ver}(-\log \dD)\rightarrow T_\xX(-\log \dD) \rightarrow p_2^* T_{\PP^n}\rightarrow 0
\end{equation*} 
\begin{equation*}\label{eqn:verticaly}
0 \rightarrow T_\yY^{\ver}(-\log \dD)\rightarrow T_\yY(-\log \dD) \rightarrow f^*p_2^*T_{\PP^n}\rightarrow 0.
\end{equation*}

We conclude this subsection by recalling the following relationship between the log sheaves defined above.

\begin{proposition}[cf. Propositions 3.3 and 3.5 in \cite{CRY}]\label{prop:log-sheaves-relationship}
Let $\xX,\dD$ and $\yY$ be as above. Then, for a general $b\in B$  
\begin{equation*}
\iota_b^*N_{\yY/\xX}(\log \dD) \simeq N_{f_b/\PP^n}(\log D_b),
\end{equation*}
and there is the following short exact sequence of sheaves on $Y_b$:
\begin{equation*}\label{eqn:log-sheaves-relationship}
0 \rightarrow \iota_b^*T_\yY^{\ver}(-\log \dD) \rightarrow \iota_b^*f^*T_\xX^{\ver}(-\log \dD) \rightarrow N_{f_b/\PP^n}(\log D_b) \rightarrow 0.
\end{equation*}
\end{proposition}

Throughout this paper we use the following setup. 

\begin{setup}\label{setup}
{\em 
Let $D\subset \PP^n$ be a divisor consisting of a general pair of hypersurfaces of degrees $d_1$ and $d_2$, and let $f:Y \rightarrow \PP^n$ be a nonconstant map from a smooth projective curve $Y$ of geometric genus $g=g(Y)$ such that 
$f(Y)\nsubset D$, $\left|f^{-1}(D)\right|=i$, and $\deg f^*H=e$ for some $i,e.$ Take $\xX,\dD$ and $\yY$ to be the families over $B$ constructed in \S\ref{sec:hilb-scheme}. 
}
\end{setup}

\subsection{Lazarsfeld kernel bundles}\label{sec:kernel-bundles}
Since we are interested in the sheaf $N_{f_b/\PP^n}(\log D_b)$, we will now study the vertical log tangent sheaf $T_\xX^{\ver}(-\log \dD)$ by comparing it to variants of the Lazarsfeld kernel bundles.

\begin{definition}[cf. \cite{Lazbundle}]
\label{def:Lazarsfeld}
{ 
The {\em Lazarsfeld kernel bundle} $M_d$
on $\PP^n$ is the kernel of the section evaluation map $H^0(\PP^n,\OO_{\PP^n}(d))\otimes \OO_{\PP^n}\rightarrow \OO_{\PP^n}(d)$.
}
\end{definition}

The main goal of this subsection is to prove the following proposition, which allows us to view the Lazarsfeld kernel bundle $M_1$ as a ``simple building block'' for the log normal sheaf $N_{f_b/\PP^n}(\log D_b)$ if the curve $f_b:Y_b\rightarrow \PP^n$ possibly violates the algebraic hyperbolicity of $(\PP^n,D_b)$. 

\begin{proposition}\label{prop:map-from-m}
As in Setup~\ref{setup},
there is a natural generically 
surjective map 
\[ \iota_b^*f^*p_2^*M_1^{\oplus s}\rightarrow N_{f_b/\PP^n}(\log D_b)\]
at a general $b\in B.$
\end{proposition}

Before proving Proposition~\ref{prop:map-from-m}, we prove some preliminary results. The following two results show how one can obtain the map in Proposition~\ref{prop:map-from-m} -- by considering a natural embedding of a direct sum of Lazarsfeld kernel bundles into a large subbundle of $T_\xX^{\ver}(-\log \dD)$ (Proposition~\ref{prop:lazarsfeld-bundles}), and composing it with a natural surjection from a direct sum of $M_1$'s induced by multiplication of polynomials (Lemma~\ref{lem:m1-to-md}).

\begin{lemma}\label{lem:m1-to-md}
For some $s>0$, there is a surjection 
\begin{equation}
\label{eqn:m1-to-md}
M_1^{\oplus s}\rightarrow \bigoplus\limits_{i=1}^2 M_{d_i}
\end{equation}
given by multiplication by polynomials.
\end{lemma}

\begin{proof}
As in \cite[Proposition 3.8]{CRY}, any pair $(P_1,P_2)$ of polynomials in $\bigoplus\limits_{i=1}^2H^0(\PP^n,\OO_{\PP^n}(d_i-1))$ induces a map $M_1\rightarrow \bigoplus\limits_{i=1}^2 M_{d_i}$ that is locally given by $s \mapsto (s\cdot P_1,s\cdot P_2).$ 
Taking a direct sum of such maps induced by generic polynomials in $\bigoplus\limits_{i=1}^2H^0(\PP^n,\OO_{\PP^n}(d_i-1))$ gives us a surjective map $M_1^{\oplus s}\rightarrow \bigoplus\limits_{i=1}^2 M_{d_i}$. 
\end{proof}

\begin{proposition}\label{prop:lazarsfeld-bundles}
There is a sheaf $M_{d_1,d_2}$ on $\xX$ such that
\begin{equation}
\label{eqn:inclusions-lazbundles}
\bigoplus\limits_{i=1}^2 p_2^*M_{d_i} \subset M_{d_1,d_2} \subset T_\xX^{\ver}(-\log \dD)
\end{equation}
and both inclusions induce trivial quotients.
\end{proposition}

\begin{proof}
There is a short exact sequence 
\begin{equation}
    0 \rightarrow T_\xX^{\ver}(-\log \dD) \rightarrow T_B \rightarrow N'_{\dD/\xX} \rightarrow 0,
\end{equation}
where $T_B\simeq \bigoplus\limits_{i=1}^2 H^0(\PP^n,\OO_{\PP^n}(d_i)) \otimes \OO_\xX.$ 
This short exact sequence fits into the following generically 
exact, commuting diagram, explained below.
\begin{equation}\label{diag:lazarsfeld-1}
\xymatrix{
& & & &0 \ar[d]\\
& &0 \ar[d] & &p_2^*\OO_{\PP^n} \ar[d]^{\mu_1}\\
&0 \ar[r] &M_{d_1,d_2} \ar[d] \ar[r] &\bigoplus\limits_{i=1}^2 H^0(\PP^n,\OO_{\PP^n}(d_i)) \otimes \OO_\xX \ar[r]^{\quad \qquad \mu_2} \ar[d]^{=} &p_2^*\OO_{\PP^n}(d_1+d_2) \ar[r] \ar[d] &0\\
&0 \ar[r] &T_\xX^{\ver}(-\log \dD) \ar[d] \ar[r] &\bigoplus\limits_{i=1}^2 H^0(\PP^n,\OO_{\PP^n}(d_i))\otimes \OO_\xX \ar[r] &N'_{\dD/\xX} \ar[r] \ar[d] &0\\
& &\kK_1 \ar[d] & &0\\
& &0
}
\end{equation}
Locally at $(p,b)\in \xX,$ the map $\mu_1$ is given by multiplication by $G_1G_2$ where $G_1G_2=0$ is the defining equation for the pair of smooth hypersurfaces corresponding to $b$. 
The 
surjective map $\mu_2$ is locally given by \linebreak $(s_1,s_2)\mapsto s_1G_2+s_2G_1$, and we denote its kernel by $M_{d_1,d_2}.$ 
By the snake lemma, the quotient $\kK_1$ in the left column is isomorphic to $p_2^*\OO_{\PP^n}\simeq \OO_\xX.$

The top row in \eqref{diag:lazarsfeld-1} fits into another 
exact, commuting diagram, explained below.
\begin{equation}\label{diag:lazarsfeld-2}
\xymatrix{
& & & &0 \ar[d]\\
& &0 \ar[d] & &p_2^*\OO_{\PP^n} \ar[d]^{\mu_3}\\
&0 \ar[r] &\bigoplus\limits_{i=1}^2 p_2^*M_{d_i}\ar[d] \ar[r] &\bigoplus\limits_{i=1}^2 H^0(\PP^n,\OO_{\PP^n}(d_i))\otimes \OO_\xX\ar[r] \ar[d]^{=} &\bigoplus\limits_{i=1}^2 p_2^*\OO_{\PP^n}(d_i) \ar[r] \ar[d]^{\mu_4} &0\\
&0 \ar[r] &M_{d_1,d_2} \ar[d] \ar[r] &\bigoplus\limits_{i=1}^2 H^0(\PP^n,\OO_{\PP^n}(d_i))\otimes \OO_\xX \ar[r]^{\quad \qquad\mu_2} &p_2^*\OO_{\PP^n}(d_1+d_2) \ar[r] \ar[d] &0\\
& &\kK_2 \ar[d] & &0\\
& &0
}
\end{equation}
The top row is obtained from the direct sum of the defining exact sequences for each $M_{d_i}$. 
The map $\mu_3$ is locally given by multiplication by $(G_1,-G_2)$. 
The 
surjective map $\mu_4$ is locally given by $(f_1,f_2)\mapsto f_1G_2+f_2G_1$. 
Again by the snake lemma, the quotient $\kK_2$ in the left column is isomorphic to $p_2^*\OO_{\PP^n}\simeq \OO_\xX.$
\end{proof}

\begin{proof}[Proof of Proposition~\ref{prop:map-from-m}]
By Proposition~\ref{prop:lazarsfeld-bundles}, there is an inclusion of $\bigoplus\limits_{i=1}^2p_2^*M_{d_i}$ into $T_\xX^{\ver}(-\log \dD)$.
Consider the restriction 
\[ \alpha:\bigoplus\limits_{i=1}^2f^*p_2^*\,M_{d_i}\rightarrow N_{f/\xX}(\log \dD)\] 
of the surjective map $f^*T_\xX^{\ver}(-\log \dD) \rightarrow N_{f/\xX}(\log \dD)$ to $\bigoplus\limits_{i=1}^2f^*p_2^*M_{d_i}$, and its pullback
\begin{equation}
\label{eqn:alpha_b}
\alpha_b\colon \bigoplus\limits_{i=1}^2\iota_b^*f^*p_2^*\,M_{d_i}\rightarrow N_{f_b/\PP^n}(\log D_b)
\end{equation}
at a general $b\in B$.
Let $G_1G_2=0$ be the defining equation for the pair of hypersurfaces corresponding to $b$, where $(G_1,G_2) \in \bigoplus\limits_{i=1}^2 H^0(\PP^n, \OO_{\PP^n}(d_i))$. 
Since $(G_1,G_2)$ and $(G_1,-G_2)$ are in $\iota_b^*T_{\yY}^{\ver}(-\log \dD)$, it follows from  diagrams \eqref{diag:lazarsfeld-1} and \eqref{diag:lazarsfeld-2} in the proof of Proposition~\ref{prop:lazarsfeld-bundles} that $\alpha_b$ is generically surjective.
Therefore, the composition
\[ \iota_b^*f^*p_2^*\,M_1^{\oplus s}\rightarrow \bigoplus\limits_{i=1}^2 \iota_b^*f^*p_2^*\,M_{d_i}\rightarrow N_{f_b/\PP^n}(\log D_b)\]
of $\alpha_b$ with the pullback of the surjective map \eqref{eqn:m1-to-md} is a generically surjective map. \qedhere

\end{proof}

\section{Proof of Theorem~\ref{introthm}}
\label{sec:3}

We use Proposition~\ref{prop:map-from-m} to show that if the curves in $\yY$ violate the inequalities in Theorem~\ref{introthm}, then we may identify a special family of lines in $\PP^n$ passing through general points on such curves. 
We call these lines \textit{associated lines} (see Definition~\ref{def:ass_line}) following \cite{ClemensRan2004, CRY}.
As in \cite{CRY}, we argue the existence and restrictive geometry of the associated lines for the $n\geq3$ and $n=2$ cases separately, in \S\ref{sec:assoc-lines-n>=3} and \S\ref{sec:assoc-lines-n=2}, respectively.
Finally, we exploit this property to show that the curves of any degree $e$ and geometric genus $g$ that violate algebraic hyperbolicity for a very general hypersurface pair $D$ form a bounded family (Propositions~\ref{prop:bounded-family-pn} and \ref{prop:bounded-family-p2}). In particular all such curves must be contained in a proper locus in $\PP^n$ (\S\ref{sec:proof-main-theorem}).

We continue following the notation of Setup~\ref{setup} from \S\ref{sec:hilb-scheme} throughout.

\subsection{Associated lines in the \texorpdfstring{$n\geq3$}{n >= 3} case}
\label{sec:assoc-lines-n>=3}

The main goal of this section is to prove the following proposition, which states that curves violating algebraic hyperbolicity form a bounded family when $n\geq 3$. 

\begin{proposition}\label{prop:bounded-family-pn}
Let $n\geq 3$, and let $D\subset \PP^n$ be a reduced divisor consisting of a very general pair of hypersurfaces of degrees $d_1$ and $d_2$ such that $d_1+d_2=2n.$
Then, the maps $f:Y \rightarrow \PP^n$ 
satisfying 
\begin{equation}
\label{eqn:alghyp-n>=3}
2g-2+i< e
\end{equation}
form a bounded family.
\end{proposition}

We first argue the existence of associated lines (see Definition \ref{def:ass_line}) using a similar argument as in \cite{ClemensRan2004, CRY}. 

\begin{lemma}\label{lem:assoc-line-pn}

If at a general point $(p,b)\in \yY$, the curve $f_b:Y_b\rightarrow X_b=\PP^n$ satisfies 
$$2g-2+i< e,$$
then there exists a unique line $\ell:=\ell(p,b)\subseteq X_b=\PP^n$ through $p$ whose ideal $\bigoplus\limits_{i=1}^2H^0(\PP^n,\OO_{\PP^n}(d_i)\otimes I_\ell)^{\oplus2}$ is contained in $\iota_b^*T_{\yY}^{vert}(-\log \dD)\vert_p$. 
\end{lemma}

\begin{definition}\label{def:ass_line}
We call a line $\ell(p,b)$ as in Lemma~\ref{lem:assoc-line-pn} the \emph{associated line} for $(p,b)$, when it exists. 
\end{definition}

\begin{proof}[Proof of Lemma~\ref{lem:assoc-line-pn}]
Recall from the proof of Proposition~\ref{prop:map-from-m} that there are natural maps \[f_b^*p_2^*M_1\rightarrow N_{f_b/\PP^n}(\log D_b).\]
Consider such a map 
\[m_{P_1,P_2}:f_b^*p_2^*M_1\rightarrow N_{f_b/\PP^n}(\log D_b)\] 
that is induced by a generic pair $(P_1,P_2)$ of polynomials in $\bigoplus\limits_{i=1}^2H^0(\PP^n,\OO_{\PP^n}(d_i-1))$ as in \eqref{eqn:m1-to-md}.
Since we assume that $2g-2+i<e$, this map must have rank one. Otherwise, we would have \[\deg N_{f_b/\PP^n}(\log D_b)\geq -(n-2)\deg f_b^*p_2^*M_1 = -(n-2)e,\] which implies that $2g-2+i\geq e$ by Proposition~\ref{prop:degree-log-normal}.
Since rank $N_{f_b/\PP^n}(\log D_b)=n-1\geq2$ and Proposition~\ref{prop:map-from-m} says that the map \eqref{eqn:alpha_b} is generically surjective, we can apply Lemma 2.2 in \cite{Clemens03} to conclude that the kernel $K$ of $m_{P_1,P_2}$ is independent of the polynomials $(P_1,P_2).$ 
Therefore, at a general point $(p,b)\in \yY$, the image of $K$ in $f_b^*T_\xX^{\ver}(-\log \dD)\vert_p$, which is contained in $\iota_b^*T_{\yY}^{vert}(-\log \dD)\vert_p$, is the ideal $\bigoplus\limits_{i=1}^2H^0(\PP^n,\OO_{\PP^n}(d_i)\otimes I_\ell)^{\oplus2}$ of some line $\ell\subset \PP^n$ through $p$. 
The uniqueness of such a line follows from the same dimension count argument in Lemma 3.13 in \cite{CRY}.
\end{proof}

As in \cite{CRY}, we define $\overline{U}^{\sym}_{0,d}$ to be the quotient of the moduli space $\overline{M}_{0,d+1}$ of stable rational curves with $d+1$ marked points by the $S_d$-action permuting the markings on the first $d$ marked points. 
In the setting of Lemma~\ref{lem:assoc-line-pn}, we may define a rational map 
\begin{equation}\label{eqn:theta}
\theta: \yY \dashrightarrow \overline{U}^{\sym}_{0,2n}
\end{equation}
that maps $(p,b)$ in $\yY$ to the modulus of $2n+1$ marked points on the associated line $\ell(p,b)$, where the first $2n$ unordered marked points are $D_b\cap \ell(p,b)$ and the $(2n+1)$st marked point is $p$. 
The rational map $\theta$ is undefined where $\ell(p,b)$ is a flex line to $D_b$, or where $\ell(p,b)$ is tangent to $D_b$ at $p$.

\begin{proof}[Proof of Proposition~\ref{prop:bounded-family-pn}]
As in the above discussion, due to Lemma~\ref{lem:assoc-line-pn}, we have the rational map $\theta$ as in \eqref{eqn:theta}.
It suffices to show that $\theta$ has constant image.

Let $\yY_p\subset \yY$ denote the fiber of $p_2\circ f: \yY \rightarrow \PP^n$ over a general point $p$ in $\PP^n$. 
Since $\yY$ is $GL(n+1)$-invariant, the associated line map $\yY_p \dashrightarrow \GG(1,n)$ maps onto the locus of lines passing through $p$. 
Let $\yY_{p,\ell}\subset \yY_p$ be the fiber over a general line $\ell.$
When restricted to $\yY_{p,\ell}$, we can factor $\theta$ as $\theta_1: \yY_{p,\ell} \rightarrow H^0(\ell,\OO_{\ell}(2n))$ composed with $\theta_2: H^0(\ell,\OO_{\ell}(2n))\dashrightarrow \overline{U}^{\sym}_{0,2n}.$
We will show that $\theta$ restricted to $\yY_{p,\ell}$ has zero-dimensional image. 

When $n\geq3$, it follows from Voisin's linear algebra argument (see  \cite[Lemma 3.17]{CRY} and  \cite[Lemma 3]{VoisinCorrection}) that $\bigoplus\limits _{i=1}^2H^0(\OO_{\PP^n}(d_i)\otimes I_\ell)$ is contained in the relative tangent space $T_{\theta_1}$ to $\theta_1$.
Let $G_1,G_2$ denote the polynomials defining the pair of hypersurfaces at $b$, then $(G_1,-G_2)$ is contained in $T_{\theta_1}$ as well.
Away from the undefined locus, the fibers of $\theta_2$ have dimension three, coming from the automorphisms of the line, with respect to the $GL(n+1)$-action, fixing a point. 
Note that $\dim \yY_{p,\ell} = \dim \yY - (2n-1) = \dim B -(2n-2).$
Hence, 
\begin{align*}
    \dim \theta(\yY_{p,\ell}) & \leq (\dim B - (2n-2)) - \sum\limits_{i=1}^2 h^0(\PP^n,\OO_{\PP^n}(d_i)\otimes I_\ell) - 1 - 3 \\
    &= \sum\limits_{i=1}^2 (h^0(\PP^n,\OO_{\PP^n}(d_i)) - h^0(\PP^n,\OO_{\PP^n}(d_i)\otimes I_\ell)) - (2n+2)\\
    &= \sum\limits_{i=1}^2 (d_i+1) - (2n+2) = 0. \qedhere
\end{align*}
\end{proof}

\subsection{Stability and associated lines in the \texorpdfstring{$n=2$}{n = 2} case}
\label{sec:assoc-lines-n=2}

We prove the analogous result for the $n=2$ case in this subsection.
Recall that in the inequalities \eqref{eqn:alghyp-n>=3} and \eqref{eqn:alghyp-n=2} the coefficients of $e$ are different in the $n\geq 3$ and $n=2$ cases, respectively.
In the $n=2$ case, the coefficient $\frac{1}{2}$ appears in the argument used in Lemma \ref{lem:assoc-line-p2} to exhibit associated lines. 

\begin{proposition}\label{prop:bounded-family-p2}
Let $D\subset \PP^2$ be a quartic plane curve consisting of a very general pair of curves. 
Then, the maps $f:Y \rightarrow \PP^2$ 
satisfying 
\begin{equation}
\label{eqn:alghyp-n=2}
2g-2+i< \frac{1}{2}e
\end{equation}
form a bounded family.
\end{proposition}

Before giving the proof of the above proposition, we prove some preliminary lemmas. 
First, we have to argue the existence of associated lines using a different proof than in Lemma~\ref{lem:assoc-line-pn}, since the argument there does not apply when $\rank N_{f_b/\PP^n}(\log D_b)=1.$

\begin{lemma}\label{lem:assoc-line-p2}
Let $D\subset \PP^2$ be a quartic plane curve consisting of a general pair of curves.
If at a general point $(p,b)\in \yY$, the curve $f_b:Y_b\rightarrow X_b=\PP^2$ satisfies 
\begin{equation}\label{eqn:1/2}
    2g-2+i< \frac{1}{2}e,
\end{equation} 
then there exists 
a unique associated line $\ell:=\ell(p,b)\subseteq X_b=\PP^2$ for $(p,b)$.
\end{lemma}

\begin{proof}
As in the proof of Lemma~\ref{lem:assoc-line-pn}, consider the map
\begin{equation*}
	m_{P_1,P_2}: f_b^*p_2^*M_1 \rightarrow \bigoplus\limits_{i=1}^2f_b^*p_2^*M_{d_i}\xrightarrow{\alpha_b} N_{f_b/\PP^2}(\log D_b)
\end{equation*}
induced by a generic pair of polynomials $(P_1, P_2)\in \bigoplus\limits_{i=1}^2H^0(\PP^2,\OO_{\PP^2}(d_i-1))$, which is generically surjective since rank $N_{f_b/\PP^2}(\log D_b)=1.$ 
The existence and uniqueness of the associated line at $(p,b)$ then follows from the same reasoning in the $n=2$ case of the proof of Lemma 3.13 in \cite{CRY}.
We will explain the existence part below.
The uniqueness of the line then follows from a dimension count.

Let us denote the torsion-free part of the image of $m_{P_1,P_2}$ by $Q_{P_1,P_2}$. 
Since $f_b^*p_2^*M_1$ has rank two, both $\ker m_{P_1,P_2}$ and $Q_{P_1,P_2}$ have rank one. 
If $\deg \ker m_{P_1,P_2} \leq \deg Q_{P_1,P_2}$, then we have $2g-2+i\geq \frac{1}{2}e$, which contradicts \eqref{eqn:1/2}.
Otherwise, $0 \subset \ker m_{P_1,P_2} \subset Q_{P_1,P_2}$ is the Harder-Narasimhan filtration, so $\ker m_{P_1,P_2}$ is independent of the choice of polynomials $P_1,P_2$, and we denote it by $K.$
Therefore, locally at $p$, the image of the map
\[ f_b^*p_2^*K\xrightarrow{-\otimes \bigoplus\limits_{i=1}^2 H^0(\PP^2, \OO_{\PP^2}(d_i-1))} \bigoplus\limits_{i=1}^2f_b^*p_2^*M_{d_i}, \]
which is the ideal $\bigoplus\limits_{i=1}^2H^0(\PP^2,\OO_{\PP^2}(d_i)\otimes I_\ell)$ of some line $\ell$ in $\PP^2$ through $p$, is contained in $\iota_b^*T_{\yY}^{vert}(-\log \dD)\vert_p$. 
\end{proof}

The following lemma shows that a preliminary dimension count as in Proposition~\ref{prop:bounded-family-pn} allows us to conclude that the image of the map $\theta$ from \eqref{eqn:theta} is not large.
However, we need to prove that the image of $\theta$ is in fact constant, which we achieve by identifying the space $\iota_b^*T_{\yY}^{\ver}(-\log \dD)|_p$  of log tangent directions to $\yY$ at a general point $(p,b)$.
We do so by using results from \cite{marangone} on rank-two vector bundles on $\PP^2$ (Lemma~\ref{lem:jumping-lines}).
 
\begin{lemma}\label{lem:modulus1-p2}
In the setting of Lemma~\ref{lem:assoc-line-p2}, the image of the map $\theta$ in \eqref{eqn:theta} has dimension at most one.
\end{lemma}

\begin{proof}
We use the same notations as in the proof of Proposition~\ref{prop:bounded-family-pn}.
By Lemma~\ref{lem:assoc-line-p2}, at least a codimension one subspace of $\bigoplus\limits_{i=1}^2H^0(\PP^2,\OO_{\PP^2}(d_i)\otimes I_\ell)$ is contained in the relative tangent space $T_{\theta_1}$ to $\theta_1$. 
Note that since $\rank N_{f_b/\PP^2}(\log D_b)=1$ here, we are not able to apply Voisin's linear algebra argument as in Proposition~\ref{prop:bounded-family-pn} to conclude immediately that $\bigoplus\limits_{i=1}^2H^0(\PP^2,\OO_{\PP^2}(d_i)\otimes I_\ell)$ is contained in $T_{\theta_1}$.
Similar to Proposition~\ref{prop:bounded-family-pn}, the fibers of $\theta_2$ have dimension three away from the undefined locus, and
$\dim \yY_{p,\ell} = \dim \yY - 3 = \dim B -2.$
Hence, 
\begin{align*}
    \dim \theta(\yY_{p,\ell}) & \leq (\dim B - 2) - \left(\sum\limits_{i=1}^2 h^0(\PP^2,\OO_{\PP^2}(d_i)\otimes I_\ell) - 1\right) - 1 - 3 \\
    &= \sum\limits_{i=1}^2 (h^0(\PP^2,\OO_{\PP^2}(d_i)) - h^0(\PP^2,\OO_{\PP^2}(d_i)\otimes I_\ell)) - 5\\
    &= 2\cdot 3 - 5 = 1.\qedhere
\end{align*}
\end{proof}

Let $\EE$ be a rank-two vector bundle on $\PP^2$. 
For a line $\ell\subset \PP^2$, we have $\EE\vert_\ell\simeq \OO_{\PP^1}(a)\oplus\OO_{\PP^1}(b)$ for some integers $a$ and $b$.
By the Grauert--M\"ulich theorem, if $\EE$ is (slope) semistable, then $\EE\vert_\ell$ has a \textit{balanced splitting type} for a general line $\ell\subset \PP^2$, meaning $a$ and $b$ satisfy $\left|a-b\right|\leq 1$ (see Theorem 2.1.4 in \cite{OSS80}). 
The lines $\ell\subset \PP^2$ such that $\EE\vert_\ell$ does not have a balanced splitting type are called the \textit{jumping lines} for $\EE$.

Now consider the log tangent sheaf $T_{\PP^2}(-\log D)$ with respect to a quartic plane curve $D\subset \PP^2$ that consists of a general pair of irreducible components. 
Since $T_{\PP^2}(-\log D)$ is semistable by \cite{Gue16}, we may apply the Grauert--M\"ulich theorem to conclude that it has at most a one-dimensional family of jumping lines.
Moreover, the following lemma states that $T_{\PP^2}(-\log D)$ actually has only finitely many jumping lines. 
Since the proof is an explicit computation, see~\cite[Lemma 3.15]{CRY}, it is omitted here.

\begin{lemma}\label{lem:finite-jumping-lines}
Let $D\subset \PP^2$ be a quartic plane curve that consists of a general pair of irreducible components. 
Then, the restriction $T_{\PP^2}(-\log D)\vert_\ell$ of the log tangent sheaf to $D$ splits as $\OO\oplus \OO(-1)$ for a general line $\ell\subset \PP^2$, and as $\OO(1)\oplus \OO(-2)$ for only finitely many lines $\ell\subset \PP^2$. 
In particular, $T_{\PP^2}(-\log D)$ has only finitely many jumping lines.
\end{lemma}

By Lemma~\ref{lem:finite-jumping-lines}, we may assume that the associated lines $\ell(p,b)$ from Lemma~\ref{lem:assoc-line-p2} are not jumping lines for $T_{\PP^2}(-\log D_b)$.
In the next Lemma~\ref{lem:jumping-lines}, we identify the space $\iota_b^*T_{\yY}^{\ver}(-\log \dD)|_p$ in this case.
We will use the following notation in the proof of this lemma. 
If $(p,b)$ is a general point in $\yY$ and $G_1,G_2$ are the polynomials defining the two components of the quartic plane curve $D_b$ at $b$, then we will denote by $J$ the ideal in $\iota_b^*\ker\alpha|_p$ generated by 
\[(s\, G_{1,x}, s \, G_{2,x}),
    (s \, G_{1,y}, s \, G_{2,y}),
    (s \, G_{1,z}, s \, G_{2,z})\]
for all $s\in M_1\vert_p$.
Here, $G_{1,x}$ denotes the partial derivative of $G_1$ with respect to $x$, etc.

\begin{lemma}\label{lem:jumping-lines}
Suppose we are in the setting of Lemma~\ref{lem:assoc-line-p2}.
If at a general point $(p,b)\in \yY$, the associated line $\ell=\ell(p,b)$ is not a jumping line for $T_{\PP^2}(-\log D_b)$, then the space 
$\iota_b^*T_{\yY}^{\ver}(-\log \dD)|_p$ is generated by $\bigoplus\limits_{i=1}^2H^0(\PP^2,\OO_{\PP^2}(d_i)\otimes I_\ell)$, $J$, $(G_1,-G_2)$, and $(G_1,G_2)$. 
\end{lemma}

\begin{proof}
Let us denote by $G=G_1G_2$ the polynomial defining the quartic plane curve $D_b$ consisting of two components defined by $G_1$ and $G_2$.
From the proof of Proposition~\ref{prop:map-from-m}, we have that the map
\[\alpha_b\colon \bigoplus\limits_{i=1}^2f_b^*p_2^*M_{d_i}\rightarrow N_{f_b/\PP^2}(\log D_b)\]
is generically surjective. 
It follows from Proposition~\ref{prop:lazarsfeld-bundles} that $\iota_b^*\ker \alpha\vert_p$ has codimension three in $f_b^*T_{\xX}^{\ver}(-\log \dD)|_p$. 
Let us denote by $c$ the codimension of $\bigoplus\limits_{i=1}^2H^0(\PP^2,\OO_{\PP^2}(d_i)\otimes I_\ell) + J $ in $\bigoplus\limits_{i=1}^2H^0(\PP^2,\OO_{\PP^2}(d_i))$, which satisfies $c\geq 3.$
We will show that $c=3.$


By \cite[Theorem 10]{marangone}, the non-jumping lines for a rank two vector bundle over $\PP^2$ are exactly the (weak) Lefschetz elements for 
$\widehat{J(G)}/ J(G)$, where $J(G)$ is the Jacobian ideal of $G$ and $\widehat {J(G)}$ denotes the saturation of the ideal $J(G)$ in $\CC[x,y,z]$ with respect to the maximal ideal (see also \cite{Sernesi2013, DimcaSticlaru}). 
We will denote $A:=\widehat{J(G)}/ J(G)$ for brevity.
In particular, $\ell$ being a (weak) Lefschetz element for $A$ means that the multiplication map 
\begin{equation*}\label{eqn:mult-map-surj}
    \times \ell: A_3 \rightarrow A_4    
\end{equation*}
is surjective, where $A_k$ denotes the degree-$k$ part of $A$.
Without loss of generality, let $p=V(y,z)$ and $\ell=V(z)$.
The surjectivity of this multiplication map means that 
$\widehat{J(G)}_4$ is generated by $\widehat{J(G)}_3\cdot z$ and $xG_x,$ $xG_y,$ $xG_z,$ $yG_x,$ $yG_y,$ $yG_z.$
It follows that $H^0(\PP^2,\OO_{\PP^2}(4))$ is generated by $H^0(\PP^2,\OO_{\PP^2}(4)\otimes I_\ell)$ and $xG_x,$ $xG_y,$ $xG_z,$ $yG_x,$ $yG_y,$ $yG_z$, and that $\bigoplus\limits_{i=1}^2H^0(\PP^2,\OO_{\PP^2}(d_i))$ is generated by $\bigoplus\limits_{i=1}^2H^0(\PP^2,\OO_{\PP^2}(d_i)\otimes I_\ell)$ and 
\begin{multline}
    \label{eqn-list}
    (G_1,-G_2), (x\, G_{1,x}, x \, G_{2,x}),
    (x \, G_{1,y}, x \, G_{2,y}),
    (x \, G_{1,z}, x \, G_{2,z}),\\
    (y\, G_{1,x}, y \, G_{2,x}),
    (y \, G_{1,y}, y \, G_{2,y}),
    (y \, G_{1,z}, y \, G_{2,z}).
\end{multline}
Since $\dim \bigoplus\limits_{i=1}^2H^0(\ell,\OO_\ell(d_i))=6$, the elements in \eqref{eqn-list} admit a relation modulo $\bigoplus\limits_{i=1}^2H^0(\PP^2,\OO_{\PP^2}(d_i)\otimes I_\ell)$ that involves at least one of 
$$(G_1,-G_2), (x\, G_{1,x}, x \, G_{2,x}),
    (x \, G_{1,y}, x \, G_{2,y}),
    (x \, G_{1,z}, x \, G_{2,z}).$$ 
Therefore, the codimension is $c=3$, and $\iota_b^*\ker \alpha|_p$ is generated by $\bigoplus\limits_{i=1}^2H^0(\PP^2,\OO_{\PP^2}(d_i)\otimes I_\ell)$ and $J$. 
It then follows that $\iota_b^*T_{\yY}^{\ver}(-\log \dD)|_p$ is generated by $\bigoplus\limits_{i=1}^2H^0(\PP^2,\OO_{\PP^2}(d_i)\otimes I_\ell)$, $J$, $(G_1,-G_2)$, and $(G_1,G_2)$. 
\end{proof}

We now prove the main result of this subsection.

\begin{proof}[Proof of Proposition~\ref{prop:bounded-family-p2}]
By Lemma~\ref{lem:modulus1-p2}, the image of $\theta$ has dimension at most one. 
We will prove that $\theta$ is in fact constant.
Suppose that the image is one-dimensional.
As before, let $\yY_p\subset \yY$ denote the fiber of $p_2\circ f: \yY \rightarrow \PP^2$ over a general point $p$ in $\PP^2$. 
Then the restriction of $\theta$ to $\yY_p$ factors as $$\yY_p \dashrightarrow B\times \GG(1,2)_p\dashrightarrow \overline{U}^{\sym}_{0,4},$$ where $\GG(1,2)_p$ denotes the space of lines passing through $p$. 
Let $Z$ denote the preimage of $\theta(\yY_p)$ in $B\times \GG(1,2)_p.$ 
Since we assume that $\theta(\yY_p)$ has dimension one, $Z$ has codimension one in $B\times \GG(1,2)_p$. 
So, the map $Z\rightarrow B$ is generically finite and induces an isomorphism of tangent spaces at a general point. 
By Lemma~\ref{lem:finite-jumping-lines}, we may assume that the general associated line $\ell$ is not a jumping line for $T_{\PP^2}(-\log D_b)$. 
Since the relative tangent space to the map $Z\rightarrow \overline{U}^{\sym}_{0,4}$ is exactly $\bigoplus\limits_{i=1}^2H^0(\PP^2,\OO_{\PP^2}(d_i)\otimes I_\ell) + J + (G_1,-G_2) + (G_1,G_2)$, which by Lemma~\ref{lem:jumping-lines} is equal to $T_\yY^{\ver}(-\log \dD)\mid_{(p,b)}$, $\yY$ must be a fiber of the map to $\overline{U}^{\sym}_{0,4}.$
\end{proof}

\subsection{Proof of Theorem~\ref{introthm}}
\label{sec:proof-main-theorem}

Finally, we prove that the main theorem follows from Propositions~\ref{prop:bounded-family-pn} and \ref{prop:bounded-family-p2}.

\begin{proof}[Proof of Theorem~\ref{introthm}]

By Proposition~\ref{prop:bounded-family-pn} (resp. Proposition~\ref{prop:bounded-family-p2}), the maps $f:Y\rightarrow \PP^n$ for $n\geq 3$ satisfying $2g-2+i<e$ (resp. the maps $f:Y\rightarrow \PP^2$ satisfying $2g-2+i<\frac{1}{2}e$) form a bounded family.
Suppose for the sake of contradiction that this family of curves dominates $\PP^n$, then $N_{f/\PP^2}(\log D_b)$ would be generically globally generated. 
By Proposition~\ref{prop:degree-log-normal}, this would imply that $2g-2+i\geq e$ , which is a contradiction.
\end{proof}

\bibliographystyle{alpha}
\bibliography{references}

\end{document}